\documentclass[10pt,openany]{article} 
\usepackage{amsmath,theorem,amsfonts,amssymb,amscd,bbm, eulervm} 
\usepackage[english]{babel}
\usepackage[verbose,includemp=false,paperwidth=15.50cm,paperheight=21cm,body={10.50cm,16cm},rmargin=1.5cm,twoside]{geometry} 
\usepackage[english]{babel} 
\usepackage[latin1]{inputenc} 
\usepackage{bm}

\theoremstyle{change}
\usepackage{hyperref}
\hypersetup{
    colorlinks = true,
    linkcolor = blue,
    anchorcolor = red,
   citecolor = blue,
    filecolor = red,
    pagecolor = red,
    urlcolor = blue}

\theoremstyle{change}
\theorembodyfont{\rmfamily}
\newtheorem{claim}{}[section]

\def\ttU{{\tt U}}

\def\bmtx{\begin{matrix}}
\def\emtx{\end{matrix}}

\def\bff{{\mathbf f}}

\def\sgn{\mathrm{sgn}}

\def\ovf{{\overline{f}}}

\def\bff{{\bf f}}

\def\QQ{\mathbb Q}
\def\1{{\mathbbm 1}}

\def\cocoa{{\hbox{\rm C\kern-.13em o\kern-.07em C\kern-.13em o\kern-.15em A}}}

\def\ttp{{\tt p}}

\def\ovD{\overline{D}}

\def\End{\mathrm{End}}

\def\w2M{\bigwedge^2M}

\def\wM{\bigwedge M}

\def\w{\wedge }
\def\bw{\bigwedge }

\protect
\protect
\protect
\protect

\protect
\protect

\def\sra{\rightarrow}
\def\lra{\longrightarrow}

\def\proof{\noindent{\bf Proof.}\,\,}
\def\qed{{\hfill\vrule height4pt width4pt depth0pt}\medskip}
\def\be{\begin{equation}}
\def\ee{\end{equation}}
\def\bclm{\begin{claim}}
\def\eclm{\end{claim}}
\def\beqn{\begin{eqnarray}}
\def\eeqn{\end{eqnarray}}
\def\beqn*{\begin{eqnarray*}}
\def\eeqn*{\end{eqnarray*}}

\title{Remarks on the Cayley-Hamilton theorem}

\author {L.~Gatto \and I.~Scherbak}

\date{}

\begin{document}

\maketitle

\abstract{\scriptsize{ \noindent We revisit the classical theorem by Cayley and Hamilton,
``{\em each endomorphism is a root of its own characteristic polynomial}'', 
  from the point of view of {\em Hasse--Schmidt derivations on an exterior algebra}.}}

\section {Formulation of the Main Result} \label{S2}

Let $A$ be  a commutative ring with unit,  $M$  a free $A$-module of rank $r$ and 
$\wM=\bigoplus_{j=0}^r\bw^jM$ its exterior algebra. 

To any endomorphism $f$ of $M$ we associate (see Section~\ref{S4}) the unique $A$-module homomorphism 
\[ 
\bff(t):\wM\lra \wM[[t]],\ \ \bff(t)=\sum_{j\geq 0}f_jt^j, \ \ f_j\in\End_A(\wM),
\] 
such that
\begin{enumerate}
\item[$\bullet$] $\bff(t)(\alpha\w \beta)=\bff(t)\alpha\w\bff(t)\beta,  \ \  \forall \alpha,\beta\in\wM;$ 

\item[$\bullet$]  $\bff(t)_{|_M}=\1_M+\sum_{j\geq 1}f^jt^j,$
\end{enumerate}
where $\1_M$ is the identity endomorphism of $M$.

Let 
\be
\det(\1_Mt-f)=t^r-e_1t^{r-1}+\ldots+(-1)^re_r\label{ee}
\ee
be the characteristic polynomial of $f$.
We prove, in Section~\ref{S4}, the following generalization of the Catyley-Hamilton theorem.
\vspace{-10pt}
\bclm{\bf Theorem.}\label{thmCH} {\em 
For each $j=1,\ldots, r$, the sequence $\{f_k\}_{k\geq 0}$ satisfies the following {\em linear recurrent relation of order} $j$,
\vspace{-3pt}
\[
{f_{i+j}}_{|_{\bw^{r-j+1}\hskip-3pt M}}\hskip-3pt -e_1{f_{i+j-1}}_{|_{\bw^{r-j+1}\hskip-3pt M}}\hskip-2pt+\hskip-2pt\,\ldots+(-1)^je_j{f_i}_{|_{\bw^{r-j+1}\hskip-3pt M}}\hskip-2pt=0,
 \]
 for all $i\geq 0$.}
\eclm 

Let $E_r(t):=\det\left(\1_M-ft\right)$. Using (\ref{ee}), we get explicitly
\be E_r(t)=1-e_1t+\ldots+(-1)^re_rt^r.\label{Ert}\ee

 For each $i\geq 0$, define
\be \ttU_i(\bff)=f_i-e_1f_{i-1}+\ldots+(-1)^re_rf_{i-r}\in \End_A(\wM),\label{ttU}\ee
with the convention that $f_k=0$ if $k<0$.

\bclm{\bf Corollary.} {\em We have
\be
\bff(t)={\ttU_0(\bff)+\ttU_1(\bff)t+\ldots+\ttU_{r-1}(\bff)t^{r-1}\over E_r(t)}\,.\label{eq:for5}
\ee
Moreover,  for each $j=1,\ldots, r$ and}  $\forall\alpha\in \bw^{r-j+1}M$, 
\[
\bff(t)(\alpha)={\ttU_0(\bff)(\alpha)+\ttU_1(\bff)(\alpha)t+\ldots+\ttU_{j-1}(\bff)(\alpha)t^{j-1}\over E_r(t)}.
\]
\eclm
\bclm{\bf Remark.}  
In the case $j=r$,  Theorem~\ref{thmCH} gives 
\[
{f_{i+j}}_{|_M}-e_1{f_{i+j-1}}_{|_M}+\ldots+(-1)^je_j{f_i}_{|_M}=0,\ \forall i\geq 0.
 \]
Thus the sequence $\{f^k\}_{k\geq 0}$, where $f^0:=\1_M$, satisfies the linear recurrence  relation of order $r$,
\[
f^{i+r}-e_1f^{i+r-1}+\ldots+(-1)^re_rf^i=0,\ \ \forall i\geq 0.
\]
This is the classic Cayley-Hamilton theorem. 
\eclm

\claim{}If $A$ contains the rational numbers, then the {\em formal Laplace transform} $L$   sending series $\sum_{n\geq 0}a_nt^n$ to $\sum_{n\geq 0}a_nn!t^n$, is invertible, \cite{gatlak}. In this case we can make some additional observations. 

Consider  the  series:
\[
u_{-j}(t):=L^{-1}\left({t^j\over E_r(t)}\right)\in A[[t]], \quad 0\leq j\leq r-1,
\]
which, as it was shown in~\cite{GS1}, form an $A$-basis of solutions to the linear ODE
\be
\left\{\begin{matrix} &y^{(r)}(t)-e_1y^{(r-1)}(t)+\ldots+(-1)^re_ry(t)=0,\cr\cr
&\hskip-126pt y(t)\in A[[t]].\end{matrix}\right.\label{ODE}
\ee
Equality  (\ref{eq:for5}) then implies
\be\sum_{n\geq 0} f_n{t^n\over n!}=\sum_{j=0}^{r-1}\ttU_j(\bff)u_{-j}(t).
\label{sol}\ee
\bclm{\bf Corollary.} {\em If $A$ is a $\QQ$--algebra,
then $\sum_{n\geq 0} f_n\displaystyle{t^n\over n!}$ solves {\em (\ref{ODE})} in
$\End_A(\wM)[[t]]$.} \qed
\eclm
  By restricting (\ref{sol}) to $M$ and recalling that $\displaystyle{\sum_{n\geq 0} f^n{t^n\over n!}=\exp(ft)}$, one obtains the refinement by Leonard \cite[1996]{Leon1} and Liz \cite[1998]{Liz} of Putzer's method \cite[1966]{Pu} to compute the exponential of a complex valued square matrix.
  \bclm{\bf Corollary.}\label{exp} {\em If $A$ is a $\QQ$--algebra and $f\in\End_A(M)$ has the characteristic polynomial {\em (\ref{ee})}, then
\[
\exp(ft)=v_0(t)\1_M+v_1(t)f+\cdots+v_{r-1}(t)f^{r-1},
\]
where $v_j(t)$  is the unique solution to {\em (\ref{ODE})} in $A[[t]]$ with the initial condition $v_j^{(i)}(t)=\delta_{ij}$, $0\leq i,j\leq r-1$}.
\eclm
\section{Derivations on  Exterior Algebra}\label{S3}
\claim{\bf Preliminaries.} \label{Pre} (See \cite{G1,SCGA})  Fix an $A$-basis   $b_0,\ldots,b_{r-1}$  of  the ring $M$. Consider $\wM=\bigoplus_{j=0}^r\bw^jM$, where $\bw^0M=A$,  and $\bw^jM$,  $1\leq j\leq r$, is the  $A$-module generated by $\{b_{i_1}\wedge\ldots\wedge b_{i_j}\}$ with the relation: 
\[
b_{i_{\sigma(1)}}\w\ldots\w b_{i_{\sigma(j)}}=\sgn(\sigma)b_{i_1}\w\ldots\w b_{i_j},\quad \sigma\in S_j,
\]
where $S_j$ denotes the permutation group on $j$ elements.
In particular, $\bw^1M=M$.  The exterior algebra structure on $\wM$ is given by  juxtaposition 
\[
\left\{\begin{matrix}\wedge&:&\bw^hM\times \bw^kM&\lra& \bw^{h+k}M,\cr\cr
&&(b_H,b_K)&\longmapsto&b_H\w b_K,\end{matrix}\right.
\]
where $b_H:=b_{i_1}\w\ldots\w b_{i_h}$ and $b_K:=b_{j_1}\w\ldots\w b_{j_k}$ are monomials in $\bw^hM$ and $\bw^kM$ respectively.

Let $t$ be an indeterminate over $\wM$ and denote by  $\wM[[t]]$ and $\End_A(\wM)[[t]]$ the corresponding  rings of formal power series with coefficients in $\wM$ and
$\End_A(\wM)$ respectively. If 
$D(t)=\sum_{i\geq 0}D_it^i$,  $\tilde{D}(t)=\sum_{j\geq 0}\tilde{D}_jt^j\in \End_A(\wM)[[t]]$, their product
 is defined as:
\[
D(t)\tilde{D}(t)\alpha=D(t)\sum_{j\geq 0}\tilde{D}_j\alpha\cdot t^j=\sum_{j\geq 0}(D(t)\tilde{D}_j\alpha)\cdot t^j,\ \ \forall\alpha\in\wM.
\]

Given $D(t)\in \End_A(\wM)[[t]]$, we use the same notation for the induced  $A$-homomorphism  $D(t):\wM\sra\wM[[t]]$ mapping  $\alpha\in\wM$ to $D(t)\alpha=\sum_{i\geq 0}D_i\alpha\cdot t^i\in \wM[[t]]$. 

The formal power series $D(t)=\sum_{i\geq 0}D_it^i$  is {\em invertible} in  $\End_A(\wM)[[t]]$ (i.e. there exists $\ovD(t)\in End_A(\wM)[[t]]$ such that $D(t)\ovD(t)=\ovD(t)D(t)=\1_{\wM}$),  if and only if $D_0$ is an automorphism of $\wM$.
If $D(t)$ is invertible, we shall write its inverse as
$
\ovD(t)=\sum_{i\geq 0}(-1)^i\ovD_i t^i$. 
With this convention, $D(t)\ovD(t)=\ovD(t)D(t)=\1_{\wM}$ if and only if
\be
\ovD_0D_j-\ovD_1D_{j-1}+\ldots+(-1)^j\ovD_jD_==0, \ \forall j\geq 1. \label{eq:dinvdbar}
\ee

\bclm{\bf Proposition}\label{equivl} {\em The following statements are equivalent:}
\begin{enumerate}
\item[{ i)}] $D(t)(\alpha\w \beta)=D(t)\alpha\w D(t)\beta$, $\ \forall \alpha,\beta\in \wM$;
\item[{ ii)}]  $D_i(\alpha\w \beta)=\sum_{j=0}^iD_{j}\alpha\w D_{i-j}\beta$, $\ \forall i\geq 0$.
\end{enumerate}
\eclm
\proof
i)$\Rightarrow$ ii). By definition of $D(t)$,  write i) as
\be
\sum_{i\geq 0}D_i(\alpha\w\beta)t^i=\sum_{j_1\geq 0}D_{j_1}\alpha\cdot t^{j_1}\w\sum_{j_2\geq 0}D_{j_2}\beta \cdot t^{j_2}.\label{eq:CH02}
\ee
Hence $D_i(\alpha\w \beta)$ is the coefficient of $t^i$ on the right hand side of~(\ref{eq:CH02}), which is  $\sum_{j_1+j_2=i}D_{j_1}\alpha\w D_{j_2}\beta=\sum_{j=0}^iD_{j}\alpha\w D_{i-j}\beta$.

\medskip

ii)$\Rightarrow$ i)  We have
\begin{eqnarray}
D(t)(\alpha\w \beta)&=&\sum_{i\geq 0}D_i(\alpha\w \beta)t^i=\sum_{i\geq 0}(\sum_{i_1+i_2=j}D_{i_1}\alpha\w D_{i_2}\beta)t^i\cr&=& \sum_{i_1\geq 0}D_{i_1}\alpha\cdot t^{i_1}\w \sum_{i_2\geq 0}D_{i_2}\beta\cdot t^{i_2}\cr&=&D(t)\alpha\w D(t)\beta. \hskip90pt \qed
\end{eqnarray}
\bclm{\bf Definition.} (Cf.~\cite{G1}) {\em  Let $D(t)\in \End_A(\wM)[[t]]$. The induced map $D(t): \wM\sra\wM[[t]]$  is called a  {\em Hasse--Schmidt derivation on}  $\wM$ {\em (}$HS$-{\em derivation} for short{\em )}, if it satisfies  the  {\em (}equivalent{\em )} conditions of Proposition~\ref{equivl}}. 
\eclm
We denote by $HS(\wM)$ the set of all HS-derivation on $\wM$.

\claim{\bf Remark.} If $D(t)\in HS(\wM)$,  then $D_1$ is an $A$--derivation of $\wM$, i.e. the usual Leibniz's rule $
D_1(\alpha\w\beta)=D_1\alpha\w \beta+\alpha\w D_1\beta.
$ holds.

\bclm{\bf Proposition.} (Cf.~\cite{G1, {SCGA}})  {\em The product of two $HS$-derivations is a $HS$-derivation. The inverse of a $HS$-derivation is a $HS$-derivation.}
\eclm
 \proof  For the product of $HS$-derivations  $D(t)$ and $\tilde{D}(t)$, the statement i) of  Proposition~\ref{equivl} holds:
 \begin{eqnarray*}
 D(t)\tilde{D}(t)(\alpha\w\beta)\hskip-6pt&=&\hskip-6pt D(t)(\sum_{j\geq0}\sum_{j_1+j_2=j}\tilde{D}_{j_1}\alpha\w\tilde{D}_{j_2}\beta)t^j\\
&=&\hskip-6pt\sum_{j\geq 0}\sum_{j_1+j_2=j}D(t) D_{j_1}\alpha\cdot t^{j_1}\w D(t)D_{j_2}\beta\cdot t^{j_2}\\
&=&\hskip-6pt D(t)\tilde{D}(t)\alpha\w D(t)\tilde{D}(t)\beta.\hskip66pt  
 \end{eqnarray*}
Similarly, for $\ovD(t)$, the inverse of the derivation $ D(t)$, we have
 \begin{eqnarray*}
\hskip24pt \ovD(t)(\alpha\w \beta)&=&\ovD(t)(D(t)\ovD(t)\alpha\w D(t)\ovD(t)\beta)\cr
&=&(\ovD(t)D(t))(\ovD(t)\alpha\w\ovD(t)\beta)\\
&=&\ovD(t)\alpha\w\ovD(t)\beta.\hskip 84pt\qed
 \end{eqnarray*}

The following property of $HS$-derivations on the exterior algebra was called  in \cite{SCGA}  the {\em integration by parts} formula.
\bclm{\bf Proposition.}  {\em  If  $\ovD(t)$ is the inverse of a derivation $D(t)$, then }
\be
 D(t)\alpha\w\beta=D(t)\alpha\w D(t)\ovD(t)\beta=D(t)(\alpha\w\ovD(t)\beta).  \label{eq:intprt}
\ee
\qed
\eclm
\claim{\bf Remark.} Equating the coefficients of $t^j$ on the left and right hand sides of~(\ref{eq:intprt}) gives an equivalent set of conditions: for any $j\geq 0$,
\be
D_j\alpha\w \beta=D_j(\alpha\w\beta)-D_{j-1}(\alpha\w\ovD_1\beta)+\ldots+(-1)^j\alpha\w\ovD_j\beta.
\ee
\vskip-90pt
\bclm{\bf Proposition.} \label{indder} {\em For any $g(t)=\sum_{j\geq 0}g_jt^j:M\sra \wM[[t]]$ there exists a unique $HS$-derivation $D(g;t)$ on $\wM$ such that}\linebreak
$D(g;t)_{|_M}=g(t).$
\eclm
\proof
For the chosen $A$-basis of  the ring $M$  (see the beginning of Section~\ref{Pre}) we necessarily have $D(g;t)(b_j)=g(t)b_j$, $0\leq j\leq r-1$. Hence, by Proposition~\ref{equivl} i), the action of $D(g;t)$ is uniquely defined on all the basis vectors of $\wM$:
\be
D(g;t)(b_{i_1}\w\ldots\w b_{i_j})=g(t)b_{i_1}\w\ldots\w g(t)b_{i_j},\ \ 1\leq j\leq r.
\ee
\qed

\section{Proof of the Main Results}\label{S4}

Let $f\in \End_A(M)$. Denote by $\overline{\bff}(t)=\sum_{j\geq0}(-1)^j\ovf_jt^j$  the unique $HS$-derivation of $\wM$ (see Proposition~\ref{indder}) extending 
\[\left(\1_M-ft\right):M\sra M[t]\subseteq \wM[[t]].\]
\vskip-90pt
\bclm{\bf Proposition.} \label{vanish} {\em  For each $i=1,\ldots,r$, the restriction  of
$\ovf_j$ to $\bw^iM$ vanishes, for  all $j>i$. }
\eclm
\proof Induction on $i$. For $i=1$ the statement holds, by definition of $\overline{\bff}(t)$.  Assume the statement true for $1\leq k\leq i-1$, and consider 
$\alpha\in\bw^iM$.  Due to $A$-linearity, its enough to take  $\alpha=m_1\w\ldots\w m_i$,
where $m_1,\ldots, m_i\in M$. For $j>i$ we have, according to Proposition~\ref{equivl} ii),
\[
\ovf_j(m_1\w\ldots\w m_{i})=\sum_{k=0}^j\ovf_k(m_1)\w\ovf_{j-k}(m_2\w\ldots\w m_{i-1}).\]
But each summand on the right hand side vanishes. Indeed, if $k=0,1$, then $j-k>i-1$ and by the induction hypothesis $\ovf_{j-k}(m_2\w\ldots\w m_{i-1})=0$. If $k\geq 2$, then $\ovf_k(m_1)=0$.  \qed

In particular, $\ovf_j$ vanishes on the entire $\wM$, for $j>r$.
\bclm{\bf Corollary.} \ \  $\overline{\bff}(t)=\1_{\wM}-\ovf_1t+\ldots+(-1)^r\ovf_rt^r$. \qed 
\eclm
Let $\bff(t)=\sum_{j\geq 0}f_jt^j$ be the inverse $HS$-derivation of $\overline{\bff}(t)$.
\bclm{\bf Proposition.} {\em For each $i\geq 1$, }
$f_i(m)=f^i(m), \forall m\in M$.
\eclm
\proof  Induction on $i$. Notice that $f_1(m)=\ovf_1(m)=f(m)$. Now assume  $f_k(m)=f^k(m)$, for $1\leq k\leq i-1$. According  to (\ref{eq:dinvdbar}), we have
\[
f_i(m)=\ovf_1f_{i-1}(m)-\ovf_2f_{i-2}(m)+\ldots+(-1)^{i-1}\ovf_i(m).
\]
By the induction hypothesis, $f_k(m)\in M$ for $1\leq k\leq i-1$, therefore  Proposition~\ref{vanish} implies the vanishing of all terms on the right hand side, but the first one. Then
\[
 f_i(m)=\ovf_1f_{i-1}(m)=\ovf_1(f^{i-1}(m))=f(f^{i-1}(m))=f^i(m). \ \ \ \qed
\]

Take  the basis element  $\zeta:=b_0\w b_1\w\ldots\w b_{r-1}$ of $\bw^rM$
(see Section~\ref{Pre}). It is unique up to the multiplication by an invertible in $A$. It is a common eigenvector of all of $f_i$ and $\ovf_j$. 
Recall that $e_1,\ldots,e_r\in A$ stay for the coefficients of the characteristic polynomial of $f$ and of $E_r(t)$, see (\ref{ee}), (\ref{Ert}).  It turns out that they are eigenvalues of ${\ovf_i}_{|_{\bw^rM}}$.
\bclm{\bf Proposition.} \label{Ee} {\em We have} $\ovf_i(\zeta)=e_i\zeta$, $1\leq i\leq r$.
\eclm
\proof Recall that $\overline{\bff}(t)$ is the $HS$-derivation extending $\1_M-ft$. Hence, like in the proof of Proposition~\ref{indder},
 \begin{eqnarray*}
\overline{\bff}(t)(\zeta)&=&\overline{\bff}(t)(b_0\wedge \ldots\wedge b_{r-1})\\
&=&\left(\1_M-ft\right)b_0\wedge \ldots\wedge\left(\1_M-ft\right)b_{r-1}\\
&=&\det\left(\1_M-ft\right)(b_0\wedge \ldots\wedge b_{r-1})=\det\left(\1_M-ft\right)(\zeta).
\end{eqnarray*}
Thus  the eigenvalue of $\ovf_i$ on $\zeta$ is the coefficient of $t^i$ in 
\[\det\left(\1_M-ft\right)=E_r(t),\] 
see   (\ref{Ert}).   \qed

\medskip
Define endomorphisms $\ttU_i(\bff)$, $i\geq 0$, via the equality
\be
E_r(t)\bff(t)=\sum_{i\geq 0}\ttU_i(\bff)t^i.\label{eq:predefui}
\ee
Comparing the coefficients of $t^i$ on the both sides gives the explicit formula
 (\ref{ttU}).
\bclm{\bf Lemma.}\label{LH} {\em  We have $\ttU_j(\bff)\zeta=0$ for all $j>0$.}
\eclm
\proof
Denote by $h_j$  the eigenvalue of $f_j$ on $\zeta$, and write 
\[ f_i(\zeta)=h_i\zeta,\ \ \bff(t)\zeta=H_r(t)\zeta, \ \ H_r(t)=1+\sum_{j\geq 1}h_jt^j.\] 
By construction, $E_r(t)H_r(t)=1$. Equivalently, we get
\[h_j-e_1h_{j-1}+\ldots+(-1)^je_j=0,\ \ j\geq 1. \]
Therefore
\begin{eqnarray*}
\ttU_j(\bff)\zeta&=&(f_j-\sum_{i=1}^r(-1)^ie_jf_{j-i})\zeta\\
&=&(h_j-e_1h_{j-1}+\ldots+(-1)^je_j)\zeta=0. \ \qed
\end{eqnarray*}

\bclm{\bf Lemma.}\label{lemma} {\em We have}
\be
\ttU_i(\bff)\alpha\w\beta=\sum_{j=0}^i(-1)^j\ttU_{i-j}(\bff)(\alpha\w \ovf_j(\beta)).\label{eq:eq9}
\ee
\eclm
\proof According to definition~(\ref{eq:predefui}), we write
\be
\sum_{i\geq 0}(\ttU_i(\bff)\alpha)t^i\w \beta=E_r(t)\bff(t)\alpha\w\beta,\label{eq:lftsd}
\ee
then apply integration by parts formula~(\ref{eq:intprt}), and again use  ~(\ref{eq:predefui}),
\begin{eqnarray}
&=&E_r(t)\bff(t)(\alpha\w\overline{\bff}(t)\beta))=\sum_{i\geq 0}\ttU_i(\bff)(\alpha\w\overline{\bff}(t)\beta)t^j\cr
&=&\sum_{i\geq 0}\ttU_i(\bff)(\alpha\w\sum_{j\geq 0}(-1)^j\ovf_j(\beta) )t^{j+i}\cr
&=&\sum_{i\geq 0}\left(\sum_{j= 0}^i\ttU_{i-j}(\bff)(\alpha\w (-1)^j\ovf_j(\beta))\right)t^{i}.
\label{eq:equalUr}
\end{eqnarray}
Comparing the coefficients of $t^i$ on the left hand side of~(\ref{eq:lftsd}) and  the right hand side of~(\ref{eq:equalUr}) gives ~(\ref{eq:eq9}).\qed

\claim{\bf Proof of Theorem~\ref{thmCH}.}  We shall show that $\ttU_k(\bff)\alpha=0$,  for any $k\geq j$  and any $\alpha\in \bw^{r-j+1}M$. 
For such $\alpha$ and for any $\beta\in\bw^{j-1}M$, we have $\alpha\w \beta\in \bw^rM$. Write $\ttU_k(\bff)\alpha\w\beta$, according to Lemma~\ref{lemma}, as
\[\ttU_k(\bff)(\alpha\w\beta)-\ttU_{k-1}(\bff)(\alpha\w \ovf_1(\beta))+\ldots+(-1)^k\alpha\w \ovf_k(\beta).\]
We see that each term but the very last is zero, by Lemma~\ref{LH}. The very last term also is zero, by Proposition~\ref{vanish} as $k>j-1$.
 Since $\ttU_k(\bff)\alpha\w\beta=0$ for any $\beta\in \bw^jM$, then $\ttU_k(\bff)\alpha=0$, for all $\alpha\in \bw^{r-j+1}M$. Hence the statement holds. \qed

In particular $\ttU_j(\bff)$ vanishes on the entire exterior algebra $\wM$ for all $j\geq r$.   

\medskip
Now we restrict each $\ttU_k(\bff)$, $0\leq k\leq r-1$, to $M$ getting: 
\[\ttU_k(\bff)_{|_M}=\ttp_k(f)=f^k-e_1f^{k-1}+\ldots+(-1)^ke_k.\]
Then (\ref{eq:predefui})  takes the form,
\be
\sum_{j\geq 0}f^jt^j={\1_M+\ttp_1(f)t+\ldots+(-1)^{r-1}\ttp_{r-1}(f)t^{r-1}\over E_r(t)}. \label{eq:prexp}
\ee
\claim{\bf Proof of Corollary~\ref{exp}.} \label{FLT} Given a  $\QQ$-algebra $A$, the {\em formal Laplace transform} $L:A[[t]]\sra A[[t]]$  and its inverse $L^{-1}$ act as follows, see~\cite{gatlak}, 
\[
L\,\sum_{n\geq 0}a_nt^n=\sum_{n\geq 0}n!a_nt^n,\ \ \
L^{-1}\sum_{n\geq 0}c_nt^n=\sum_{n\geq 0}c_n{t^n\over n!}.
\]

Consider
\[
u_{-j}:=L^{-1}\left(t^j\over E_r(t)\right), \ \ 0\leq j\leq r-1.
\]
Applying $L^{-1}$ to~(\ref{eq:prexp}),  we get  the following expression for $\exp(ft)$,
\[
\exp(ft)=u_0+\ttp_1(f)u_{-1}+\ldots+\ttp_{r-1}(f)u_{-r+1}.
\]
It will be convenient to re-write the series $u_0,u_{-1},\ldots,u_{-r+1}$ in terms of
\[ H_r(t)={1\over E_r(t)}=1+\sum_{j\geq 0}h_jt^j,\]
introduced in the proof of Lemma~\ref{LH}. We obtain
\[
u_{-j}=L^{-1}(t^jH_r(t))=\sum_{n\geq 0}h_{n-j}{t^n\over n!},\ \ 0\leq j\leq r-1.
\]
According to~\cite{GS1}, these series form an $A$-basis of solutions to the  ODE (\ref{ODE}) in $A[[t]]$. Hence $\exp(ft)$  solves this ODE  in $A[[t]]$. 

Take the {\em standard} $A$-basis of solutions, $\{v_j(t)\}_{0\leq j\leq r-1}$, where $v_j(t)$  denotes the unique solution to (\ref{ODE}) in $A[[t]]$ satisfying $v_j^{(i)}(t)=\delta_{ij}$, $0\leq i,j\leq r-1$. In  the standard basis,  the coefficients are the initial conditions of the solution. In the case of $\exp(ft)$, these are $\1_M,f,\ldots, f^{r-1}$. \qed



\begin{center}
\begin{scriptsize}
\hskip11pt \begin{tabular}{lcl}
{\sc Dip. di Scienze Matematiche}&&{\sc School of Math. Sciences}\\
{\sc Politecnico di Torino, Italy}&&{\sc Tel Aviv University, Israel}\\
 letterio.gatto@polito.it&& scherbak@post.tau.ac.il
\end{tabular}
\end{scriptsize}
\end{center}
\noindent

\end{document}